# Convergence and asymptotic normality of variational Bayesian approximations for exponential family models with missing values

**Bo Wang**
Department of Statistics
University of Glasgow
Glasgow G12 8QQ
Scotland, U.K.

**D. M. Titterington**
Department of Statistics
University of Glasgow
Glasgow G12 8QQ
Scotland, U.K.

## Abstract

We study the properties of variational Bayes approximations for exponential family models with missing values. It is shown that the iterative algorithm for obtaining the variational Bayesian estimator converges locally to the true value with probability 1 as the sample size becomes indefinitely large. Moreover, the variational posterior distribution is proved to be asymptotically normal.

## 1   INTRODUCTION

Variational Bayes approximations have recently been applied to complex models involving incomplete-data for which computational difficulties arise with the ideal Bayesian approach. Such models include hidden Markov models and mixture models; see for example Attias (1999, 2000); Beal (2003); Ghahramani and Beal (2000); Humphreys and Titterington (2000, 2001); MacKay (1997); Penny and Roberts (2000); Wang and Titterington (2004b). In these earlier contributions, the approximations were shown empirically to be convergent and effective. However little has been done to investigate their theoretical properties, and the purpose of this paper is to go some way to rectifying this.

Hall, Humphreys and Titterington (2002) initiated a discussion of these issues and proved that, for certain Markov models, the parameter estimator obtained by maximising the variational lower bound function is asymptotically consistent provided the proportion of all values that are missing tends to zero. Later we proved in Wang and Titterington (2003) that it is not always the case that a fully factorised form of variational posterior, which includes the factorisation of the joint probability function for the hidden states, provides an asymptotically consistent estimator as the 'sample size' becomes large. We demonstrated this in particular in the context of linear state space models, in which the above sufficient condition obviously does not hold. On the other hand we showed in Wang and Titterington (2004a) that variational Bayes estimators for certain mixture models are asymptotically efficient for large sample sizes.

In this paper we study the properties of variational approximation algorithms for more general models, namely exponential family models with missing values. Exponential families include cases such as Gaussian, gamma, Poisson, Dirichlet and Wishart distributions, and exponential family models with missing values contain many models of practical interest as particular cases, such as Gaussian mixtures, hidden Markov models and linear state space models. Beal (2003) and Ghahramani and Beal (2000) applied the variational Bayesian method to these models and derived the iterative algorithm for learning the approximate posterior distributions of the latent states and the model parameters. The numerical expriments therein show empirically that this algorithm is convergent and efficient. In this paper we derive the iterative procedure for obtaining the variational Bayesian estimator, we provide analytical proofs of local convergence of the procedure as the sample size tends to infinity, and we show that the variational posterior distribution for the parameters is asymptotically normal.

## 2   EXPONENTIAL FAMILY MODELS WITH MISSING VALUES AND VARIATIONAL APPROXIMATIONS

We consider the following exponential family models with missing values. Suppose that $\Theta$ is an open subset of $\mathbb{R}^m$, that $\mathcal{P} = \{P_\theta : \theta \in \Theta\}$ is a family of probability distributions on a measurable space $(\Omega, \mathcal{F})$, and that $x$ and $y$ are sampled from the natural exponential family

with density

$$p(x, y|\theta) = f(x, y) \exp\{\theta^\top u(x, y) - \psi(\theta)\}, \quad (1)$$

with $x$ taking values in $\mathbb{R}^d$ and $y$ in $\mathbb{R}^p$, where $\theta \in \Theta$ is the unknown parameter, and $\psi(\cdot) : \Theta \mapsto \mathbb{R}$ is six-times continuously differentiable and has positive definite Hessian matrix on $\Theta$. The parameter $\theta$ has a conjugate prior to the complete-data likelihood (1), with density

$$p(\theta|\alpha_0, \beta_0) = h(\alpha_0, \beta_0) \exp\{\theta^\top \beta_0 - \alpha_0 \psi(\theta)\}, \quad (2)$$

where $h$ is a normalising constant satisfying

$$h(\alpha, \beta)^{-1} = \int_\Theta \exp\{\theta^\top \beta - \alpha \psi(\theta)\} d\theta, \quad (3)$$

and $\alpha_0 \in \mathbb{R}, \beta_0 \in \mathbb{R}^m$ are the hyperparameters of the prior.

**Remark 1.** *The models of the forms* (1) *and* (2) *include most latent-variable models of practical interest. A simple example is when $x$ is sampled from a univariate Gaussian distribution with mean $\theta_1$ and variance 1 and $y = x + w$, where $w$ is sampled from another Gaussian distribution, independent of $x$, with mean $\theta_2$ and variance 1. The joint probability density is*

$$p(x, y|\theta_1, \theta_2) = \exp\{-\frac{1}{2}(y - x)^2 + \theta_1 x + \theta_2(y - x)$$
$$- \frac{1}{2}(\theta_1^2 + \theta_2^2) - \log(2\pi)\}.$$

*And the parameters $\theta_1$ and $\theta_2$ have independent Gaussian prior distributions with the same variance.*

Suppose that only $y$ is observable whereas $x$ is latent. We have a data-set consisting of a random sample of size $n$, with $Y = (y_1, y_2, \dots, y_n)$ and $X = (x_1, x_2, \dots, x_n)$. In the Bayesian framework we want to infer the posteriors over both the parameters and the hidden states. Unfortunately exact Bayesian inference is generally time-consuming, if not impossible, especially for large dimensionality $m$. Therefore approximation is usually necessary in these cases. In the variational approach, the true posterior $p(X, \theta|Y)$ is approximated by the variational distribution $q(X, \theta)$, which factorises as $q(X, \theta) = q_X(X)q_\theta(\theta)$, and is chosen to maximise the negative free energy

$$\int q(X, \theta) \log \frac{p(\theta, X, Y)}{q(X, \theta)} d\theta dX, \quad (4)$$

equivalent to minimising the Kullback-Leibler divergence between the exact and approximate distributions of $\theta$ and $X$, given $Y$.

The negative free energy (4) can be maximised using the following iterative procedure (Beal (2003);

Ghahramani and Beal (2000)). In turn, the following two stages are performed.

(i) Optimise $q_\theta(\theta)$ for fixed $\{q_{x_i}(x_i), i = 1, \dots, n\}$, defined in (ii) below. This step results in

$$q_\theta(\theta) = h(\alpha, \beta) \exp\{\theta^\top \beta - \alpha \psi(\theta)\}, \quad (5)$$

where $\alpha$ and $\beta$ are the hyperparameters of the variational posterior and are updated by

$$\alpha = n + \alpha_0, \quad \beta = \sum_{i=1}^n r_i + \beta_0, \quad \text{and } r_i = \langle u(x_i, y_i) \rangle_{x_i}.$$
$$(6)$$

Here $\langle \cdot \rangle_{x_i}$ denotes the expectation under $q_{x_i}(x_i)$.

(ii) Optimise $q_X(X)$ for fixed $q_\theta(\theta)$. This leads to the factorised form $q_X(X) = \prod_{i=1}^n q_{x_i}(x_i)$, where

$$q_{x_i}(x_i) = f(x_i, y_i)g(\theta, y_i)$$
$$\exp\{\langle \theta \rangle_\theta^\top u(x_i, y_i) - \psi(\langle \theta \rangle_\theta)\}, \quad (7)$$

in which $g(\theta, y_i)$ is a normalising constant satisfying

$$g(\theta, y_i)^{-1} = \int f(x_i, y_i)$$
$$\exp\{\langle \theta \rangle_\theta^\top u(x_i, y_i) - \psi(\langle \theta \rangle_\theta)\} dx_i, \quad (8)$$

and $\langle \cdot \rangle_\theta$ denotes the expectation under $q_\theta(\theta)$.

## 3 THE ITERATIVE ALGORITHM AND ITS CONVERGENCE

We define the variational Bayesian estimator $\hat{\theta}$ of the parameter $\theta$ as

$$\hat{\theta} = \int_\Theta \theta q_{\text{pos}}(\theta) d\theta,$$

where $q_{\text{pos}}$ is the variational posterior density of $\theta$, given by the limiting form of $q_\theta(\theta)$ that results from the above iterative procedure. For the exponential family distribution (5) the corresponding variational Bayesian estimator is

$$\hat{\theta} = \int_\Theta \theta q_\theta(\theta) d\theta = -\frac{D_\beta h(\alpha, \beta)}{h(\alpha, \beta)}.$$

(Throughout the paper, $D\Psi$ and $D^2\Psi$ denote the gradient and the Hessian of $\Psi$. When ambiguity exists, the specific variable of differentiation appears as a subscript of the symbol $D$ and $D^2$.)

Thus, the procedure in the previous section can be used to derive the following algorithm for obtaining the variational Bayesian estimate of $\theta$: starting with some initial value $\theta^{(0)}$, successive iterates are defined inductively by

$$\theta^{(k+1)} \triangleq \Phi_n(\theta^{(k)}) = -\frac{D_\beta h(\alpha,\beta)}{h(\alpha,\beta)}, \qquad (9)$$

where $\alpha$ and $\beta$ are given as in (6), and

$$\begin{aligned} q_{x_i}(x_i) = {}& f(x_i,y_i)g(\theta^{(k)},y_i) \\ & \exp\{(\theta^{(k)})^\top u(x_i,y_i) - \psi(\theta^{(k)})\}, \\ g(\theta^{(k)},y_i)^{-1} = {}& \int f(x_i,y_i) \\ & \exp\{(\theta^{(k)})^\top u(x_i,y_i) - \psi(\theta^{(k)})\} dx_i. \end{aligned}$$

It is of interest to investigate the questions of whether or not the algorithm (9) is convergent and, if so, what properties are possessed by the limiting value. The following theorem gives a partial answer.

**Theorem 1.** *With probability 1 as $n$ approaches infinity, the iterative procedure (9) converges locally to the true value $\theta^*$, i.e. (9) converges to $\theta^*$ whenever the starting value is sufficiently near to $\theta^*$.*

*Proof.* Define the norm of $\theta \in \mathbb{R}^m$ as $\|\theta\| \triangleq (\theta^\top \theta)^{1/2}$ and the norm of the real $m \times m$ matrix $A$ as $\|A\| \triangleq \sup_{\|\theta\|=1} \|A\theta\|$.

We first prove that, with probability 1 as $n$ approaches infinity, the operator $\Phi_n$ is locally contractive; that is, there exists a number $\lambda$, $0 \le \lambda < 1$, such that

$$\|\Phi_n(\bar\theta) - \Phi_n(\theta^*)\| \le \lambda\|\bar\theta - \theta^*\|, \qquad (10)$$

whenever $\bar\theta$ lies sufficiently near $\theta^*$.

Since $\bar\theta$ is near $\theta^*$ we can write

$$\Phi_n(\bar\theta) - \Phi_n(\theta^*) = D\Phi_n(\theta^*)(\bar\theta - \theta^*) + O(\|\bar\theta - \theta^*\|^2),$$

where $D\Phi_n(\theta^*)$ denotes the gradient of $\Phi_n(\theta)$ evaluated at $\theta^*$. It follows that

$$\begin{aligned} \|\Phi_n(\bar\theta) - \Phi_n(\theta^*)\| \le {}& \|D\Phi_n(\theta^*)\| \cdot \|\bar\theta - \theta^*\| \\ = {}& \sup_{\|\theta\|=1} |\theta^\top D\Phi_n(\theta^*)\theta| \cdot \|\bar\theta - \theta^*\|. \end{aligned}$$

Consequently, it is sufficient to show that $D\Phi_n(\theta^*)$ converges with probability 1 to a matrix which has norm less than 1.

Write $\beta$ and $r_i$ as $\beta(\theta)$ and $r_i(\theta)$ to indicate explicitly their dependence on $\theta$. From (9) one has

$$\begin{aligned} & D\Phi_n(\theta^*) \\ ={}& \frac{D_\beta h(\alpha,\beta)D_\beta^\top h(\alpha,\beta) - h(\alpha,\beta)D_\beta^2 h(\alpha,\beta)}{h^2(\alpha,\beta)} D\beta(\theta^*). \end{aligned}$$

Here $h$ and its derivatives are evaluated at $\theta^*$. For convenience we write $h(\alpha,\beta)^{-1}$ evaluated at $\theta^*$ as $\tilde h(\alpha,\beta)$, from which

$$D_\beta \tilde h(\alpha,\beta) = \int_\Theta \exp\{\theta^\top \beta - \alpha\psi(\theta)\}\theta d\theta,$$

$$D_\beta^2 \tilde h(\alpha,\beta) = \int_\Theta \exp\{\theta^\top \beta - \alpha\psi(\theta)\}\theta\theta^\top d\theta.$$

Let $b(\cdot) : \mathbb{R}^m \mapsto \mathbb{R}$ be a four-times continuously differentiable function of $\theta$ and write

$$a_n(\theta) = (1 + \frac{\alpha_0}{n})\psi(\theta) - \theta^\top(\frac{1}{n}\sum_{i=1}^n r_i + \frac{\beta_0}{n}), \qquad (11)$$

$$h_b = \int_\Theta b(\theta)\exp\{-na_n(\theta)\}d\theta. \qquad (12)$$

Since $\psi(\theta)$ is continuously differentiable and has positive definite Hessian matrix, it is obvious that $a_n(\theta)$ is also continuously differentiable and strictly convex in $\theta$. Thus, if we let $\hat\theta_n$ solve the equation

$$D\psi(\theta) = (\frac{1}{n}\sum_{i=1}^n r_i + \frac{\beta_0}{n})/(1 + \frac{\alpha_0}{n}), \qquad (13)$$

$\hat\theta_n$ is also the unique global minimiser of $a_n(\theta)$ on $\Theta$.

It is obvious that $D^2 a_n$ converges to $D^2\psi$ with probability 1 as $n \to \infty$. By Lemma 1 in Appendix B, letting $b(\theta)$ be 1, $\theta_i$ and $\theta_i\theta_j$ $(i,j=1,\dots,m)$ correspondingly in (23) and after a straightforward calculation, we obtain that, as $n$ tends to infinity, with probability 1,

$$\begin{aligned} & \frac{nD_\beta^{2,ij}\tilde h(\alpha,\beta)\tilde h(\alpha,\beta) - nD_\beta^i\tilde h(\alpha,\beta)D_\beta^j\tilde h(\alpha,\beta)}{\tilde h^2(\alpha,\beta)} \\ & \to \frac{1}{2}\sigma_\infty^{ij} = \frac{1}{2}[D^2\psi(\theta)]_{ij}^{-1}. \end{aligned} \qquad (14)$$

In Appendix A we prove that, as $n \to \infty$,

$$\frac{1}{n}D\beta(\theta^*) \to D^2\psi(\theta^*) - \mathbb{E}_{y_i}\{\mathbb{E}_{x_i}[\phi]\mathbb{E}_{x_i}[\phi^\top]\}, \ a.s.,$$

where 'a.s.' means 'almost surely' and $\phi$ is defined as

$$\phi = u(x_i,y_i) - D\psi(\theta^*). \qquad (15)$$

Therefore, combining (14) with the last limiting result we obtain that, with probability 1,

$$\begin{aligned} & D\Phi_n(\theta^*) \\ \to{}& \frac{1}{2}[D^2\psi(\theta^*)]^{-1}\big[D^2\psi(\theta^*) - \mathbb{E}_{y_i}\{\mathbb{E}_{x_i}[\phi]\mathbb{E}_{x_i}[\phi^\top]\}\big] \\ ={}& \frac{1}{2}I_m - \frac{1}{2}[D^2\psi(\theta^*)]^{-1}\mathbb{E}_{y_i}\{\mathbb{E}_{x_i}[\phi]\mathbb{E}_{x_i}[\phi^\top]\}, \end{aligned}$$

where $I_m$ denotes the $m \times m$ identity matrix.

Since $D^2\psi(\theta^*)$ is positive definite and symmetric and obviously $\mathbb{E}_{y_i}\{\mathbb{E}_{x_i}[\phi]\mathbb{E}_{x_i}[\phi^\top]\}$ is positive semidefinite and symmetric, $D\Phi_n(\theta^*) \leq \frac{1}{2}I_m$ as $n$ tends to infinity; that is, $D\Phi_n(\theta^*) - \frac{1}{2}I_m$ is negative semidefinite.

Next we show that

$$[D^2\psi(\theta^*)]^{-1}\mathbb{E}_{y_i}\{\mathbb{E}_{x_i}[\phi]\mathbb{E}_{x_i}[\phi^\top]\} \leq I_m.$$

Since $D^2\psi(\theta^*)$ is positive definite and symmetric, it is sufficient to prove that

$$\theta^\top \mathbb{E}_{y_i}\{\mathbb{E}_{x_i}[\phi]\mathbb{E}_{x_i}[\phi^\top]\}\theta \leq \theta^\top D^2\psi(\theta^*)\theta, \qquad (16)$$

for any $\theta \in \mathbb{R}^m$.

In fact, we have

$$\theta^\top \mathbb{E}_{y_i}\{\mathbb{E}_{x_i}[\phi]\mathbb{E}_{x_i}[\phi^\top]\}\theta$$
$$=\mathbb{E}_{y_i}\Big\{\sum_{j,k=1}^m \theta_j\theta_k\mathbb{E}_{x_i}[\phi_j]\mathbb{E}_{x_i}[\phi_k]\Big\}$$
$$=\mathbb{E}_{y_i}\Big\{\sum_{j,k=1}^m \theta_j\theta_k\mathbb{E}_{x_i}[\phi_j\phi_k]$$
$$\quad - \sum_{j,k=1}^m \theta_j\theta_k\mathbb{E}_{x_i}\big[(\phi_j - \mathbb{E}_{x_i}[\phi_j])(\phi_k - \mathbb{E}_{x_i}[\phi_k])\big]\Big\}$$
$$=\mathbb{E}_{y_i}\Big\{\sum_{j,k=1}^m \theta_j\theta_k\mathbb{E}_{x_i}[\phi_j\phi_k]$$
$$\quad - \theta^\top \mathbb{E}_{x_i}\big[(\phi - \mathbb{E}_{x_i}[\phi])(\phi - \mathbb{E}_{x_i}[\phi])^\top\big]\theta\Big\}$$
$$\leq \sum_{j,k=1}^m \theta_j\theta_k\mathbb{E}_{y_i}\{\mathbb{E}_{x_i}[\phi_j\phi_k]\}$$
$$=\theta^\top D^2\psi(\theta^*)\theta,$$

where the last equality is a consequence of (22).

Therefore, we obtain $0 \leq D\Phi_n(\theta^*) \leq \frac{1}{2}I_m$, and consequently the inequality (10) holds with $\lambda = 1/2$. Moreover, if we use Laplace's approximation (23) it is easy to deduce that $\Phi_n(\theta^*) = -D_\beta h(\alpha,\beta)/h(\alpha,\beta) \to \theta^*$ with probability 1 as $n$ tends to infinity.

Therefore, since the starting value is sufficiently near to $\theta^*$ we have

$$\|\theta^{(k+1)} - \theta^*\|$$
$$\leq \|\Phi_n(\theta^{(k)}) - \Phi_n(\theta^*)\| + \|\Phi_n(\theta^*) - \theta^*\|$$
$$\leq \lambda\|\theta^{(k)} - \theta^*\| + \|\Phi_n(\theta^*) - \theta^*\|,$$

and therefore the iterative procedure (9) converges locally to the true value $\theta^*$ with probability 1 as $n$ approaches infinity . $\qquad\square$

# 4 ASYMPTOTIC NORMALITY OF THE VARIATIONAL POSTERIOR DISTRIBUTION

There have been a large number of contributions about the asymptotic normality of posterior distributions associated with exponential families; see for instance Walker (1969), Heyde and Johnstone (1979), Chen (1985) and Bernardo and Smith (1994). Under appropriate conditions the (true) posterior density converges in distribution to a normal density. In this section, we show that the variational posterior distribution for the parameter $\theta$ obtained by the iterative procedure also has the property of asymptotic normality. This implies that the variational posterior becomes more and more concentrated around the true parameter value as the sample size grows.

Suppose the sample size $n$ is large. We have proved that the algorithm (9) is convergent, so there exists an equilibrium point denoted by $\tilde{\theta}_n$. It follows from (5) and (7) that, at $\tilde{\theta}_n$,

$$\tilde{\alpha}_n = n + \alpha_0, \quad \tilde{\beta}_n = \sum_{i=1}^n r_i + \beta_0, \quad r_i = \langle u(x_i, y_i)\rangle_{x_i},$$
$$q(x_i) = f(x_i, y_i)g(\tilde{\theta}_n, y_i)\exp\{\tilde{\theta}_n^\top u(x_i, y_i) - \psi(\tilde{\theta}_n)\}.$$

Therefore, the variational posterior density of $\theta$ at the equilibrium point is

$$q_n(\theta) = h(\tilde{\alpha}_n, \tilde{\beta}_n)\exp\{\theta^\top \tilde{\beta}_n - \tilde{\alpha}_n\psi(\theta)\}.$$

Let $\hat{\theta}_n$ maximise $\theta^\top \tilde{\beta}_n - \tilde{\alpha}_n\psi(\theta)$. Then we have

$$D\psi(\hat{\theta}_n) = \Big(\frac{1}{n}\sum_{i=1}^n r_i + \frac{\beta_0}{n}\Big)/\Big(1 + \frac{\alpha_0}{n}\Big).$$

By the same arguments as used in the previous section and noting that $\tilde{\theta}_n \to \theta^*$ with probability 1 by Theorem 1, we have that $\frac{1}{n}\sum_{i=1}^n r_i$ converges to $D\psi(\theta^*)$ almost surely. Since $D\psi$ is strictly increasing and continuous, $\hat{\theta}_n \to \theta^*$ with probability 1 as $n$ tends to infinity.

Define

$$L_n(\theta) \triangleq \log q_n(\theta) = \log h(\tilde{\alpha}_n, \tilde{\beta}_n) + \theta^\top \tilde{\beta}_n - \tilde{\alpha}_n\psi(\theta).$$

Then we have

$$\Sigma_n \triangleq -[D^2 L_n(\hat{\theta}_n)]^{-1} = [(n + \alpha_0)D^2\psi(\hat{\theta}_n)]^{-1}.$$

Denote by $B(\theta, \varepsilon)$ the open ball of radius $\varepsilon$ centred at $\theta$. According to Chen (1985), under the assumption of the consistency of $\hat{\theta}_n$ for $\theta^*$, the posterior density $q_n$ converges in distribution to $\mathcal{N}(\hat{\theta}_n, \Sigma_n)$ if the following basic conditions hold.

(C1) "Steepness". $\sigma_n^2 \to 0$ with $P_{\theta^*}$-probability 1 as $n \to \infty$, where $\sigma_n^2$ is the largest eigenvalue of $\Sigma_n$.

(C2) "Smoothness". For any $\varepsilon > 0$, there exists an integer $N$ and $\delta > 0$ such that, for any $n > N$ and $\theta \in B(\theta, \delta)$, $D^2 L_n(\theta)$ exists and satisfies

$$I_m - A(\varepsilon) \le D^2 L_n(\theta)[D^2 L_n(\hat{\theta}_n)]^{-1} \le I_m + A(\varepsilon), \ a.s.,$$

where $A(\varepsilon)$ is an $m \times m$ symmetric positive semidefinite matrix whose largest eigenvalue tends to zero with $P_{\theta^*}$-probability 1 as $\varepsilon \to 0$.

(C3) "Concentration". For any $\delta > 0$,

$$\int_{B(\theta, \delta)} q_n(\theta) d\theta \to 1$$

with $P_{\theta^*}$-probability 1 as $n$ tends to infinity.

In fact, since $\hat{\theta}_n \to \theta^*$, the components of $D^2 \psi(\hat{\theta}_n)$ are bounded above and away from 0 almost surely if $n$ is large enough, so the largest eigenvalue of $\Sigma_n$ tends to 0.

(C2) is obvious because $D^2 L_n(\theta)[D^2 L_n(\hat{\theta}_n)]^{-1} = D^2 \psi(\theta)[D^2 \psi(\hat{\theta}_n)]^{-1}$ and $\psi(\cdot)$ is continuously differentiable.

From Kass et al. (1990), assumption (iii) in Appendix B is stronger than (C3). Therefore all the conditions are verified.

## Acknowledgement

This work was supported by a grant from the UK Science and Engineering Research Council. The authors gratefully acknowledge the reviewers for their valuable comments.

## Appendix A

In the appendix we prove that the following convergences hold:

$$\frac{1}{n}D\beta(\theta^*) \to D^2\psi(\theta^*) - \mathbb{E}_{y_i}\{\mathbb{E}_{x_i}[\phi]\mathbb{E}_{x_i}[\phi^\top]\}, \ a.s., \tag{17}$$

$$\frac{1}{n}\beta(\theta^*) \to D\psi(\theta^*), \ a.s. \tag{18}$$

In fact, from (6) we have that $D\beta(\theta) = \sum_{i=1}^n Dr_i(\theta)$ and

$$
\begin{aligned}
&Dr_i(\theta^*) = \int u(x_i, y_i)D_\theta^\top q_{x_i}(x_i)dx_i \\
&= \int u(x_i, y_i)f(x_i, y_i)D_\theta^\top g(\theta^*, y_i) \\
&\qquad \cdot \exp\{\theta^{*\top}u(x_i, y_i) - \psi(\theta^*)\}dx_i \\
&+ \int u(x_i, y_i)f(x_i, y_i)g(\theta^*, y_i)\exp\{\theta^{*\top}u(x_i, y_i) - \psi(\theta^*)\} \\
&\qquad \cdot [u^\top(x_i, y_i) - D^\top\psi(\theta^*)]dx_i \\
&= \int [u(x_i, y_i) - D\psi(\theta^*)]f(x_i, y_i)D_\theta^\top g(\theta^*, y_i) \\
&\qquad \cdot \exp\{\theta^{*\top}u(x_i, y_i) - \psi(\theta^*)\}dx_i \\
&+ \int f(x_i, y_i)g(\theta^*, y_i)\exp\{\theta^{*\top}u(x_i, y_i) - \psi(\theta^*)\} \\
&\qquad \cdot [u(x_i, y_i) - D\psi(\theta^*)][u^\top(x_i, y_i) - D^\top\psi(\theta^*)]dx_i,
\end{aligned}
$$

where in the last equality we used the fact that

$$
\begin{aligned}
&\int f(x_i, y_i)D_\theta g(\theta^*, y_i)\exp\{\theta^{*\top}u(x_i, y_i) - \psi(\theta^*)\}dx_i \\
&\quad + \int [u(x_i, y_i) - D\psi(\theta^*)]f(x_i, y_i)g(\theta^*, y_i) \\
&\qquad \cdot \exp\{\theta^{*\top}u(x_i, y_i) - \psi(\theta^*)\}dx_i = 0, \tag{19}
\end{aligned}
$$

which is obtained by differentiating, with respect to $\theta$,

$$\int f(x_i, y_i)g(\theta^*, y_i)\exp\{\theta^{*\top}u(x_i, y_i) - \psi(\theta^*)\}dx_i = 1.$$

Since it follows from (8) that

$$
\begin{aligned}
&D_\theta g(\theta^*, y_i) \\
&= -\int f(x_i, y_i)\exp\{\theta^{*\top}u(x_i, y_i) - \psi(\theta^*)\} \\
&\qquad \cdot [u(x_i, y_i) - D\psi(\theta^*)]dx_i \\
&\quad \cdot \left\{\int f(x_i, y_i)\exp\{\theta^{*\top}u(x_i, y_i) - \psi(\theta^*)\}dx_i\right\}^{-2},
\end{aligned}
$$

equality (19) can be rewritten as

$$
\begin{aligned}
&\int [u(x_i, y_i) - D\psi(\theta^*)]f(x_i, y_i)g(\theta^*, y_i) \\
&\qquad \cdot \exp\{\theta^{*\top}u(x_i, y_i) - \psi(\theta^*)\}dx_i \\
&\quad \cdot \int f(x_i, y_i)\exp\{\theta^{*\top}u(x_i, y_i) - \psi(\theta^*)\}dx_i \\
&= \int f(x_i, y_i)\exp\{\theta^{*\top}u(x_i, y_i) - \psi(\theta^*)\} \\
&\qquad \cdot [u(x_i, y_i) - D\psi(\theta^*)]dx_i. \tag{20}
\end{aligned}
$$

Differentiating both sides of (20) with respect to $\theta^*$, we have

$$
\begin{aligned}
&\left\{\int [u(x_i, y_i) - D\psi(\theta^*)]f(x_i, y_i)D_\theta^\top g(\theta^*, y_i) \right. \\
&\qquad \cdot \exp\{\theta^{*\top}u(x_i, y_i) - \psi(\theta^*)\}dx_i - D^2\psi(\theta^*) \\
&\quad + \int f(x_i, y_i)g(\theta^*, y_i)\exp\{\theta^{*\top}u(x_i, y_i) - \psi(\theta^*)\} \\
&\qquad \left. \cdot [u(x_i, y_i) - D\psi(\theta^*)][u^\top(x_i, y_i) - D^\top\psi(\theta^*)]dx_i\right\} \\
&\quad \cdot \int f(x_i, y_i)\exp\{\theta^{*\top}u(x_i, y_i) - \psi(\theta^*)\}dx_i \\
&\quad + \int [u(x_i, y_i) - D\psi(\theta^*)]f(x_i, y_i)g(\theta^*, y_i) \\
&\qquad \cdot \exp\{\theta^{*\top}u(x_i, y_i) - \psi(\theta^*)\}dx_i \\
&\quad \cdot \int f(x_i, y_i)\exp\{\theta^{*\top}u(x_i, y_i) - \psi(\theta)\} \\
&\qquad \cdot [u^\top(x_i, y_i) - D^\top\psi(\theta^*)]dx_i \\
&= \int f(x_i, y_i)\exp\{\theta^{*\top}u(x_i, y_i) - \psi(\theta^*)\} \\
&\quad \cdot [u(x_i, y_i) - D\psi(\theta^*)][u^\top(x_i, y_i) - D^\top\psi(\theta^*)]dx_i \\
&\quad - \int f(x_i, y_i)\exp\{\theta^{*\top}u(x_i, y_i) - \psi(\theta^*)\}dx_i \cdot D^2\psi(\theta^*). \tag{21}
\end{aligned}
$$

We define $\phi$ as in (15). The marginal distribution of $y_i$ is $\int p(x_i, y_i|\theta^*)dx_i$, and therefore it follows from the strong law of large numbers that, with probability 1,

$$
\begin{aligned}
&\frac{1}{n}\sum_{i=1}^{n}\nabla_\theta r_i\\
&\to \int\Big\{\nabla_\theta r_i\int p(x_i,y_i|\theta^*)dx_i\Big\}dy_i\\
&= D^2\psi(\theta^*)-\int\Big\{\mathbb{E}_{x_i}[\phi]\mathbb{E}_{x_i}[\phi^\top]\\
&\qquad\cdot\int f(x_i,y_i)\exp\{\theta^{*\top}u(x_i,y_i)-\psi(\theta^*)\}dx_i\Big\}dy_i\\
&= D^2\psi(\theta^*)-\mathbb{E}_{y_i}\big\{\mathbb{E}_{x_i}[\phi]\mathbb{E}_{x_i}[\phi^\top]\big\},
\end{aligned}
$$

where we have used equality (21) and the fact that

$$
\begin{aligned}
&\int f(x_i,y_i)\exp\{\theta^{*\top}u(x_i,y_i)-\psi(\theta^*)\}\\
&\cdot\big[u(x_i,y_i)-D\psi(\theta^*)\big]\big[u^\top(x_i,y_i)-D^\top\psi(\theta^*)\big]dx_idy_i\\
&= D^2\psi(\theta^*),
\end{aligned}\tag{22}
$$

and $\mathbb{E}_{x_i}$ denotes expectation under $q_{x_i}$.

Thus, we obtain (17). Derivation of (18) is similar.

## Appendix B

In this appendix we show that under our framework the Laplace approximation is justified. The proof consists of verifying the *analytical assumptions for Laplace's method* in Kass, Tierney and Kadane (1990), which are listed here for convenience. Since $a_n$ defined in (11) is of random nature, some minor revisions are made to adapt to our settings.

Suppose that $\{a_n : n = 1, 2, \dots\}$ is a sequence of six-times continuously differentiable real functions and that $b$ is a four-times continuously differentiable function of $\theta$. The pair $(\{a_n\}, b)$ is said to satisfy the *analytical assumptions for Laplace's method* if there exist positive numbers $\varepsilon$, $M$ and $\eta$, and an integer $n_0$ such that $n > n_0$ implies the following:

(i) for all $\theta \in B(\hat\theta_n, \varepsilon)$ and all $1 \le j_1, \dots, j_d \le m$ with $0 \le d \le 6$, $|\partial_{j_1\cdots j_d}a_n(\theta)| < M$ with $P_{\theta^*}$-probability 1;

(ii) $\det(D^2a_n(\hat\theta_n)) > \eta$ with $P_{\theta^*}$-probability 1;

(iii) the integral $h_b$ defined in equation (12) exists and is finite, and, for all $\delta$ for which $0 < \delta < \varepsilon$, $B(\hat\theta_n, \delta) \subseteq \Theta$,

$$
\begin{aligned}
&\big[\det(nD^2a_n(\hat\theta_n))\big]^{1/2}\int_{\Theta-B(\hat\theta_n,\delta)}b(\theta)\\
&\qquad\cdot\exp\{-n(a_n(\hat\theta_n)-a_n(\theta))\}d\theta = O(n^{-2})
\end{aligned}
$$

with $P_{\theta^*}$-probability 1; or, more strongly,

(iii') for all $\delta$ for which $0 < \delta < \varepsilon$, $B(\hat\theta_n, \delta) \subseteq \Theta$,

$$
\limsup_{n\to\infty}\sup_\theta\{a_n(\hat\theta_n)-a_n(\theta):\theta\in\Theta-B(\hat\theta_n,\delta)\}<0
$$

with $P_{\theta^*}$-probability 1.

According to Kass et al. (1990), we have the following lemma.

**Lemma 1.** *If $(\{a_n\}, b)$ satisfy the analytical assumptions for Laplace's method then*

$$
\begin{aligned}
&\int_\Theta b(\theta)\exp\{-na_n(\theta)\}d\theta\\
&= (2\pi)^{m/2}[\det(nD^2a_n)]^{-1/2}\exp\{-na_n(\hat\theta_n)\}\\
&\quad\cdot\Big\{b(\hat\theta_n)+\frac{1}{n}\Big[\frac{1}{2}\sum_{i,j=1}^m\sigma_n^{ij}b_{ij}-\frac{1}{6}\sum_{\substack{i,j=1\\k,s=1}}^m a_n^{ijk}b_s\mu_{ijks}^4\\
&\quad+\frac{1}{72}b(\hat\theta_n)\sum_{\substack{i,j,k=1\\q,r,s=1}}^m a_n^{ijk}h_n^{qrs}\mu_{ijkqrs}^6\\
&\quad-\frac{1}{24}b(\hat\theta_n)\sum_{\substack{i,j=1\\k,s=1}}^m a_n^{ijks}\mu_{ijks}^4\Big]+O(n^{-2})\Big\},\quad a.s.,
\end{aligned}\tag{23}
$$

*where $\mu_{ijks}^4$ and $\mu_{ijkqrs}^6$ are the fourth and sixth central moments of a multivariate normal distribution having covariance matrix $(D^2a_n)^{-1}$; that is,*

$$
\begin{aligned}
\mu_{ijks}^4 =\ &\sigma_n^{ij}\sigma_n^{ks}+\sigma_n^{ik}\sigma_n^{js}+\sigma_n^{is}\sigma_n^{jk},\\
\mu_{ijkqrs}^6 =\ &\sigma_n^{ij}\sigma_n^{kq}\sigma_n^{rs}+\sigma_n^{ij}\sigma_n^{kr}\sigma_n^{qs}+\sigma_n^{ij}\sigma_n^{ks}\sigma_n^{qr}\\
&+\sigma_n^{ik}\sigma_n^{jq}\sigma_n^{rs}+\sigma_n^{ik}\sigma_n^{jr}\sigma_n^{qs}+\sigma_n^{ik}\sigma_n^{js}\sigma_n^{qr}\\
&+\sigma_n^{iq}\sigma_n^{jk}\sigma_n^{rs}+\sigma_n^{iq}\sigma_n^{jr}\sigma_n^{ks}+\sigma_n^{iq}\sigma_n^{js}\sigma_n^{kr}\\
&+\sigma_n^{ir}\sigma_n^{jk}\sigma_n^{qs}+\sigma_n^{ir}\sigma_n^{jq}\sigma_n^{ks}+\sigma_n^{ir}\sigma_n^{js}\sigma_n^{kq}\\
&+\sigma_n^{is}\sigma_n^{jk}\sigma_n^{qr}+\sigma_n^{is}\sigma_n^{jq}\sigma_n^{kr}+\sigma_n^{is}\sigma_n^{jr}\sigma_n^{kq},
\end{aligned}
$$

*where $D^2a_n$ denotes the Hessian of $a_n$, its $(i,j)$-component is written as $a_n^{ij}$ and the components of its inverse are written as $\sigma_n^{ij}$; moreover, $b_s$ and $b_{ij}$ denote the components of the first- and second-order derivatives of $b$, respectively. All derivatives are evaluated at $\hat\theta_n$.*

Now we verify the assumptions (i)-(iii).

Under our assumptions, it has been shown in (18) that $\frac{1}{n}\sum_{i=1}^n r_i \to D\psi(\theta^*)$ with probability 1, so, when $n$ large enough, $\frac{1}{n}\sum_{i=1}^n r_i$ is almost surely bounded in $B(\hat\theta_n, \varepsilon)$. Since $\psi$ is continuously differentiable (i) obviously holds.

Condition (ii) is one of our assumptions.

As $n$ tends to infinity, for any $\theta \in \Theta$, $a_n(\theta)$ converges with $P_{\theta^*}$-probability 1 to

$$
a_0(\theta)=\psi(\theta)-\theta^\top D\psi(\theta^*).
$$

Since $\hat{\theta}_n$ maximises $a_n$, we have

$$\hat{\theta}_n = (D\psi)^{-1}\big(\big(\frac{1}{n}\sum_{i=1}^{n} r_i + \frac{\beta_0}{n}\big)/(1 + \frac{\alpha_0}{n})\big),$$

so it follows that, as $n$ tends to infinity, with probability 1,

$$\hat{\theta}_n \to (D\psi)^{-1}\big(D\psi(\theta^*)\big) = \theta^*.$$

Therefore, for all $\delta$ for which $0 < \delta < \varepsilon$ and $\theta \in \Theta - B(\hat{\theta}_n, \delta)$, we have that, $\forall \varepsilon_0$ satisfying $0 < \varepsilon_0 < \delta/2$, there exists an integer $N$ such that, if $n > N$, it holds that, for all $\theta \in \Theta$,

$$|a_n(\theta) - a_0(\theta)| < \varepsilon_0,$$

$$\|\hat{\theta}_n - \theta^*\| < \varepsilon_0, \ a.s.$$

$$|a_0(\hat{\theta}_n) - a_0(\theta^*)| < \varepsilon_0, \ a.s.$$

Thus,

$$\begin{aligned}
a_n(\hat{\theta}_n) - a_n(\theta) &= a_n(\hat{\theta}_n) - a_0(\hat{\theta}_n) + a_0(\hat{\theta}_n) - a_0(\theta^*) \\
&\quad + a_0(\theta^*) - a_0(\theta) + a_0(\theta) - a_n(\theta) \\
&< a_0(\theta^*) - a_0(\theta) + 3\varepsilon_0, \ a.s.,
\end{aligned}$$

so that

$$\begin{aligned}
&\sup\{a_n(\hat{\theta}_n) - a_n(\theta) : \theta \in \Theta - B(\hat{\theta}_n, \delta)\} \\
&\leq \sup\{a_0(\theta^*) - a_0(\theta) : \theta \in \Theta - B(\hat{\theta}_n, \delta)\} + 3\varepsilon_0 \\
&\leq \sup\{a_0(\theta^*) - a_0(\theta) : \theta \in \Theta - B(\theta^*, \delta - \varepsilon_0)\} \\
&\quad + 3\varepsilon_0, \ a.s., \qquad\qquad\qquad\qquad\qquad (24)
\end{aligned}$$

since $B(\theta^*, \delta - \varepsilon_0) \subset B(\hat{\theta}_n, \delta)$.

Since $a_0(\cdot)$ is strictly convex, for $\theta \in \Theta - B(\theta^*, \delta - \varepsilon_0)$, we have $a_0(\theta) - a_0(\theta^*) > c$, where $c = \inf\{a_0(\theta) - a_0(\theta^*) : \theta$ lies in the boundary of $B(\theta^*, \delta/2)\} > 0$. Consequently, we get

$$\sup\{a_0(\theta^*) - a_0(\theta) : \theta \in \Theta - B(\theta^*, \delta - \varepsilon_0)\} \leq -c.$$

Combining the last estimate with (24) we have that, $\forall \varepsilon_0$ satisfying $0 < \varepsilon_0 < \delta$, there exists an integer $N$ such that $n > N$ implies

$$\sup\{a_n(\hat{\theta}_n) - a_n(\theta) : \theta \in \Theta - B(\hat{\theta}_n, \delta)\} \leq -c + 3\varepsilon_0, \ a.s.;$$

that is, (iii') holds.